\documentclass{amsart}
\usepackage{amssymb, mathtools, fullpage, color}
\usepackage[colorinlistoftodos]{todonotes}
\usepackage[colorlinks=true, pdfstartview=FitV,linkcolor=blue,citecolor=blue,urlcolor=blue]{hyperref}
\usepackage{bbm}
\usepackage{listings}
\usepackage[bottom]{footmisc}

\theoremstyle{definition}

\newenvironment{red}{\relax\color{red}}{\relax}
\newenvironment{blue}{\relax\color{blue}}{\hspace*{.5ex}\relax}

\newcommand{\ber}{\begin{red}}
\newcommand{\er}{\end{red}}
\newcommand{\beb}{\begin{blue}}
\newcommand{\eb}{\end{blue}}

\numberwithin{equation}{section}

\begin{document}

\title[Predicting root numbers with neural networks]{Predicting root numbers with neural networks}

\date{\today}

\author[A. Pozdnyakov]{Alexey Pozdnyakov}
\address{Department of Mathematics, University of Connecticut, Storrs, CT 06269, U.S.A.}
\email{alexey.pozdnyakov@uconn.edu}

\begin{abstract}
We report on two machine learning experiments in search of statistical relationships between Dirichlet coefficients and root numbers or analytic ranks of certain low-degree $L$-functions. The first experiment is to construct interpretable models based on murmurations, a recently discovered correlation between Dirichlet coefficients and root numbers. We show experimentally that these models achieve high accuracy by learning a combination of Mestre-Nagao type heuristics and murmurations, noting that the relative importance of these features varies with degree. The second experiment is to search for a low-complexity statistic of Dirichlet coefficients that can be used to predict root numbers in polynomial time. We give experimental evidence and provide heuristics that suggest this can not be done with standard machine learning techniques.
\end{abstract}

\maketitle
\section{Introduction}

    Due to open-source machine learning platforms such as TensorFlow \cite{TF} and publically available mathematical databases such as the $L$-functions and modular forms database \cite{lmfdb}, machine learning has become increasingly prevalent in mathematics research. In this paper, we address several questions that have emerged from previous machine learning investigations into number theory and arithmetic geometry.
    
    One of these previous investigations is the series of papers \cite{ABH, HLOa, HLOb, HLOc, HLOP, HLOPS}. In this work, the authors apply many different statistical models and machine learning techniques such as Bayesian classifiers, random forests, logistic regression, principle component analysis, and topological data analysis to various ends. One goal is to investigate open problems in number theory such as the Sato-Tate conjecture and the famous Birch and Swinnerton-Dyer (BSD) conjecture \cite{ABH, HLOa}. Another goal was to uncover novel relationships between various algebraic and analytic invariants of objects such as arithmetic curves and number fields \cite{HLOb, HLOc}. One surprising outcome of applying machine learning methods to elliptic curve data on the LMFDB was the discovery of murmurations, a correlation between the Frobenius traces $a_p(E)$ and the root numbers $w_E$ of elliptic curves. This was first reported in \cite{HLOP}, and then generalized to the setting of $L$-functions associated with higher weight modular forms or higher genus curves in \cite{HLOPS}. This has since led to multiple theorems proving the existence of and giving exact analytic expressions for such correlations \cite{BBLL, LOP, Z23}.

    Other relevant work includes \cite{KV22}, in which large convolutional neural networks (CNNs) were used to give high-accuracy predictions of the rank $r_E$ of elliptic curves from their Frobenius traces $a_p(E)$ and conductors $N_E$. These models were shown to outperform Mestre-Nagao sums, series involving $a_p(E)$ which heuristically converge to $r_E$. With this in mind, \cite{KV22} concludes with several questions including whether the CNNs have discovered new mathematics (beyond the BSD conjecture and Mestre-Nagao type heuristics), as well as whether variations of their models or other deep learning techniques can be used to classify elliptic curves more efficiently. We investigate these questions in what follows, both in the original setting of elliptic curves and in the setting of low-degree $L$-functions.
    
    We will also address a closely related question posed by Peter Sarnak \cite{S23}. Sarnak asks whether statistical correlations between $a_p(E)$ and $w_E$ can be used to predict root numbers in polynomial time with respect to $N_E$. Note that by the parity conjecture, which follows from BSD, we have that $w_E = (-1)^{r_E}$. Therefore, predicting ranks is sufficient for predicting root numbers. We also note that the standard method for computing the root number of an elliptic curve $E$ with discriminant $\Delta_E$ is to compute
    \begin{equation}
        w_E = -\prod_{p \mid \Delta_E} w_p(E),
    \end{equation}
    where $w_p(E)$ are local root numbers. This requires factoring $\Delta_E$, which prevents polynomial time computation. On the other hand, murmurations suggest that $a_p(E)$ can be used to predict $w_E$ without factoring the discriminant. It remains to check whether this can be done with a sufficiently small number of $a_p(E)$, particularly as $N_E \to \infty$. 
    
    This problem is also connected to the Möbius function $\mu(n)$, which presents an interesting machine learning challenge. As demonstrated by \cite{H21}, which attempts to use supervised learning to predict the closely related Louisville $\lambda$ function, it is difficult for standard machine learning algorithms to predict $\mu(n)$. However, one can show (e.g. \cite{HG}) that for an elliptic curve $E$ with square-free discriminant,
    \begin{equation}
        \mu(\Delta_E) = -w_E \left(\frac{-c_6}{\Delta_E}\right),
    \end{equation}
    where $c_6$ is an integer invariant that is easily computed from the Weierstrass coefficients of $E$ and $\left(\frac{\cdot}{\cdot}\right)$ denotes the Jacobi symbol. Thus a polynomial time algorithm for predicting $w_E$ will allow us to predict $\mu(n)$ in polynomial time for certain $n$. Note that \cite{BHK} proved that, conditional on the generalization Riemann hypothesis for Dirichlet $L$-functions, there is an algorithm for determining $\mu(n) = 0$ without factoring.  

    We address the questions related to \cite{KV22} by building small neural network models. These are shallow models that have the advantage of being small enough to interpret while still achieving high accuracy. Inspired by \cite{HLOPS}, we demonstrate that such models are also effective for predicting root numbers of $L$-functions of degrees 3 and 4. By inspecting the model weights, we conclude that these neural networks achieve high accuracy by learning to compute a combination of Mestre-Nagao type heuristics and murmurations, and that murmurations appear to be particularly predictive in degree 3. For questions related to \cite{S23}, we note that the grand Riemann hypothesis for certain Rankin-Selberg $L$-functions implies that $O(\log(N_E)^2)$ Frobenius traces will suffice to uniquely determine the isogeny class of $E$, and thus $w_E$. However, this procedure has high complexity, and there is no known algorithm for computing $w_E$ from anything less than $O\big(N_E^{1/2 + o(1)}\big)$ Frobenius traces. To this end, we perform an extensive hyperparameter search on large CNN models that use at most $O(\log(N_E)^3)$ Frobenius traces to predict $w_E$. From these experiments, as well as a heuristic analysis, we conclude that standard machine learning methods are unlikely to find a method for predicting root numbers in polynomial time.

\subsection*{Acknowledgements}
The author is grateful to Andrew Granville, Yang-Hui He, Kyu-Hwan Lee, Thomas Oliver, Peter Sarnak, Andrew Sutherland, and Micheal Rubinstein for helpful discussions that shaped the direction and scope of this work. They would also like to thank Andrew Granville for sharing the probabilistic heuristic in section 5, Kyu-Hwan Lee for introducing them to AI-assisted research, Peter Sarnak for suggesting the polynomial time constraint and sharing insightful notes, and Andrew Sutherland for providing computational resources, endless data, and the algorithm in section 3.
    
\section{Background}
\subsection{Arithmetic Curves \& Low Degree $L$-functions} Our primary classification problem involves elliptic curves, although this can be considered through the language of $L$-functions. In this section, we will introduce elliptic curves and their $L$-functions before giving examples of higher degree $L$-functions. For $A, B \in \mathbb{Z}$ with $4A^3 + 27B^3 \neq 0$, let
\begin{equation}
    E/\mathbb{Q} : y^2 = x^3 + Ax + B
\end{equation}
be a non-CM elliptic curve in Weierstrass form. We say that a prime $p$ is a good (resp. bad) prime for $E$ if the reduction of the Weierstrass equation mod $p$ defines a smooth (resp. singular) curve over $\mathbb{F}_p$.
The $L$-function attached to $E$ is given by the Euler product:
\begin{equation}\label{eq.ellipticL}
L(E,s)=\prod_{p \in \mathcal{P}} L_p(E,s)^{-1},
\end{equation}
in which $\mathcal{P}$ denotes the set of primes,
\begin{equation}
    L_p(E,s) = \begin{cases}1-a_p(E)p^{-s}+p^{1-2s} &\text{if $p$ is good,}\\
    1-a_p(E)p^{-s} &\text{if $p$ is bad,}\end{cases}
\end{equation}
and
\begin{equation}
    a_p(E)=p+1-\#E(\mathbb{F}_p).
\end{equation}
For a bad prime $p$, we have $a_p\in\{-1,0,1\}$ according to the reduction type \cite[Section~VII.11]{Sil1}.
It follows from the celebrated modularity theorem that $L(E,s)$ admits analytic continuation to $\mathbb{C}$ and satisfies the functional equation
\begin{equation}
    \Lambda(E,s)=w_E\Lambda(E,2-s), \label{eq:ellipticFeq}
\end{equation}
where $w_E\in\{-1,1\}$ is the root number and
\begin{equation}
    \Lambda(E,s)=N_E^{s/2}\Gamma_{\mathbb{C}}(s)L(E,s),
\end{equation}
in which $N_E$ denotes the conductor of $E$, and $\Gamma_{\mathbb{C}}(s)=\Gamma_{\mathbb{R}}(s)\Gamma_{\mathbb{R}}(s+1)=2(2\pi)^{-s}\Gamma(s)$ is the standard archimedean Euler factor attached to the complex place. Moreover, the Mordell–Weil theorem implies that the points of $E/\mathbb{Q}$ form a finitely-generated abelian group:
\begin{equation}
    E(\mathbb{Q}) \cong E(\mathbb{Q})_{\text{tors}} \times \mathbb{Z}^{r_E},
\end{equation}
where $r_E$ is the rank of $E$. Computing the rank is a difficult open problem, and the BSD conjecture proposes that one can do so via the $L$-function:
\begin{equation}
    \text{ord}_{s=1}L(E,s) = r_E.
\end{equation}
The left-hand side is called the analytic rank of $L(E,s)$, and this quantity can be considered for any $L$-function. We note that computing the analytic rank rigorously can be challenging for ranks $> 1$.\footnote{See \cite{BGZ} for an example with rank 3.} For $L$-functions of degree 2, we will be using $(a_p(E))_{p\in\mathcal{P}}$ and $N_E$ to predict $w_E$ and $r_E$. 

Next, we describe the degree 3 $L$-functions we will be using. Let $E/\mathbb{Q}$ be an elliptic curve and let $\chi$ be a primitive quadratic Dirichlet character.
For $L(E,s)$ as in equation \eqref{eq.ellipticL} and
\begin{equation}
    L(\chi,s) = \sum_{n=1}^{\infty}\chi(n)n^{-s}=\prod_{p \in \mathcal{P}}(1-\chi(p)p^{-s})^{-1},
\end{equation}
the product $L(E,s+\frac{1}{2})L(\chi,s)$ is a degree $3$ $L$-function. This $L$-function may be identified with that of an Eisenstein series on $\mathrm{GL}_3(\mathbb{R})$ \cite[Chapter~10]{Goldfeld}.
Note that the normalization of $s$ is different from that in Equation~\eqref{eq.ellipticL}, so that the product satisfies a functional equation with respect to $s\mapsto1-s$. At a good prime $p$ for $E$, the Euler polynomial for $L(\chi,s)L(E,s+\frac12)$ is
\begin{equation}
    1-\left(\chi(p)+\frac{a_p}{\sqrt{p}}\right)p^{-s}+\left(\frac{\chi(p)a_p}{\sqrt{p}}+1\right)p^{-2s}-\chi(p)p^{-3s}.
\end{equation}
At a bad prime, we get
\begin{equation}
    1-\left(\chi(p)+\frac{a_p}{\sqrt{p}}\right)p^{-s}+\frac{\chi(p)a_p}{\sqrt{p}}p^{-2s}.
\end{equation}
The conductor of such an $L$-function is $N_\chi N_E$, where $N_\chi$ is the conductor of $\chi$. The root number and analytic rank of this $L$-function come directly from $L(E, s)$. Therefore, our classification problem is to use $(\chi(p) + \frac{a_p}{\sqrt{p}})_{p \in \mathcal{P}}$ and $N_\chi N_E$ to predict $w_E$ and $r_E$.

\quad Finally, we consider degree $4$ $L$-functions arising from genus 2 curves. Let $C$ be a smooth, projective, geometrically integral curve of genus 2. For a good prime $p$ of $C$, we define the local zeta-function 
\begin{equation}
    Z(C/\mathbb{F}_p, T) = \exp\Big(\sum_{k=1}^\infty \frac{\# C(\mathbb{F}_p) T^k}{k}\Big).
\end{equation}
Since $C$ has genus 2, the local zeta-function can be written as
\begin{equation}
    Z(C/\mathbb{F}_p, T) = \frac{L_p(C, T)}{(1-T)(1-pT)},
\end{equation}
where $L_p(C,T) \in \mathbb{Z}[T]$ is given by
\begin{equation}\label{eq.Genus2Euler}
    L_p(C,T)=1+a_{1,p}(C)T+a_{2,p}(C)T^2+a_{1,p}(C)pT^3+p^2T^4.
\end{equation}
For a bad prime $p$ of $C$, we set $a_{1,p}(C) = 0$ and $a_{2,p}(C) = p$. The Hasse–Weil conjecture claims that
\begin{equation}
    L(C, s) = \prod_{p \in \mathcal{P}} L_p(C, p^{-s})^{-1} \label{eq:genus2L}
\end{equation}
admits a meromorphic continuation to $\mathbb{C}$ which satisfies a functional equation analogous to Equation~\eqref{eq:ellipticFeq}. The curve $C$ also has a conductor $N_C$ determined by its bad primes and the reduction type at those primes. In this case, we use $(a_{1,p}(C))_{p \in \mathcal{P}}$ and $N_C$ to predict the root number and analytic rank of the $L$-function in Equation~\eqref{eq:genus2L}.
\subsection{Deep learning}

In this section, we give a general introduction to the deep learning techniques we will be using. Since our goal is to predict ranks and root numbers, we will be working in the paradigm of supervised learning. We introduce some terminology to explain the general framework. We assume that we have a set of labeled data $D = \{(x_i, y_i)\}_{i=1}^N$, where $x_i \in \mathbb{R}^n$ and $y_i \in \mathbb{R}^m$ for some $n, m \in \mathbb{N}$. We choose a model $f: \mathbb{R}^n \to \mathbb{R}^m$ given by $f(x;\theta)$ where $x \in \mathbb{R}^n$ and $\theta \in \Theta$ is a choice of parameters that define the behavior of $f$. We also pick a loss function $\mathcal{L}:\Theta \to \mathbb{R}_{\geq 0}$ which is some measure of the inaccuracy of the model on the dataset $D$. For example, when $f$ is a logistic regression model, we have $\Theta = \mathbb{R}^{n+1}$ and $f(x;\theta_w, \theta_b) = \sigma(\theta_w \cdot x + \theta_b)$ where $\sigma$ is a logistic function. In this case, it is common to use categorical cross-entropy as the loss function: 
\begin{equation}
    \mathcal{L}(\theta) = - \frac{1}{N}\sum_{i=1}^N\sum_{j=1}^m \mathbbm{1}_{y_i \in C_j}  \log\left(\frac{\exp(f_j(x_i; \theta))}{\sum_{k=1}^m \exp(f_k(x_i;\theta))} \right),
\end{equation}
% For example, a common choice in regression problems is the mean square error: 
% \begin{equation}
%     \mathcal{L}(\theta) = \frac{1}{N}\sum_{i=1}^N |f(x_i;\theta) - y_i|^2.
% \end{equation}
where $f = (f_1, \ldots, f_m)$. We think of $f_j(x_i; \theta)$ as the score the model assigns to the event $x_i \in C_j$, and we think of 
\begin{equation}
   \frac{\exp(f_j(x_i; \theta))}{\sum_{k=1}^C \exp(f_k(x_i;\theta))}
\end{equation}
as the predicted probability of $x_i \in C_j$. Once we fix our dataset, model, and loss function, we train our model. Training the model refers to solving the optimizing problem
\begin{equation}\label{eq:loss_opt}
    \theta^* = \text{argmin}_{\theta \in \Theta}  \mathcal{L}(\theta).
\end{equation}
This is usually done using some version of gradient descent. For example, see \cite{KB14} for the Adam optimizer, a standard stochastic gradient descent algorithm for training deep neural networks.

Deep learning refers to the use of deep neural networks as the model $f$ in supervised learning. Standard fully connected neural networks are functions of the following form. Let $L \in \mathbb{N}$ and $N_0, \ldots, N_L \in \mathbb{N}$. Then a feed-forward neural network $f$ with activation function $\rho : \mathbb{R} \to \mathbb{R}$ is a map $f:\mathbb{R}^{N_0} \to \mathbb{R}^{N_L}$ defined by 
\begin{equation}\label{d-NN}
    f(x;\theta) = \begin{cases}
    W_1(x),  & L = 1, \\
    W_2 \circ \rho \circ W_1(x), & L =2, \\
    W_L \circ \rho \circ W_{L-1} \circ \rho \circ \ldots \circ \rho \circ W_1(x), & L \geq 3,
\end{cases}
\end{equation}
where
\begin{equation}
    W_\ell(x) = A_\ell x + b_\ell \text{ with } A_\ell \in \mathbb{R}^{N_\ell \times N_{\ell -1}} \text{  and  } b \in \mathbb{R}^{N_\ell} \text{ for all }  \ell \in \{1,\ldots, L \},
\end{equation}
and $\rho$ acts component-wise on vectors. We typically call $L$ the depth of the neural network, and we call $\max \{N_0, \ldots, N_L\}$ the width of the neural network. In this case, the parameters $\theta$ refers to the entries in $\{A_\ell\}_{\ell=1}^L$ and $\{b_\ell\}_{\ell=1}^L$. Neural networks are a general-purpose model that achieves high accuracy in problems ranging from image recognition \cite{KSH} to dynamics prediction \cite{LKB}. One reason for their success is that they are universal approximators, meaning that they can approximate any continuous real-valued function with a compact domain. We make this precise with the universal approximation theorem, proven by Cybenko in \cite{C89}. Let $C(X, Y)$ denote the set of all continuous functions $f: X \to Y$. Let $\rho \in C(\mathbb{R},\mathbb{R})$, and let $\rho \circ x$ denote $\rho$ applied to each component of $x$. Then $\rho$ is not polynomial if and only if for every $n,m \in \mathbb{N}$, compact $K \subset \mathbb{R}^n$, $f \in C(K,\mathbb{R}^m)$, and $\epsilon > 0$, there exists $k \in \mathbb{N}$, $A_1 \in \mathbb{R}^{k \times n}$, $b_1 \in \mathbb{R}^k$, and $A_2 \in \mathbb{R}^{m \times k}$ such that
\label{t-UAT}
\begin{equation}
    \sup_{x \in K} |f(x)-g(x)| < \epsilon \text{ where } g(x) = A_2  (\rho \circ (A_1  x + b_1)).
\end{equation}
Note that similar theorems exist that fix the width but require arbitrary depth. 

Another reason for the success of neural networks is their efficient training. Solving the optimization problem in Equation~\eqref{eq:loss_opt} requires a method for taking gradients of the form 
\begin{equation}
    \frac{\partial \mathcal{L}}{\partial A_\ell} \text{ and } \frac{\partial \mathcal{L}}{\partial b_\ell} \text{ } \forall \ell \in \{1, \ldots, L\},
\end{equation}
which is inefficient to do directly with neural networks of high depth. Instead, we take these gradients using the well-known backpropagation algorithm. This algorithm computes these gradients inductively using the chain rule. Let us define $a_\ell \in \mathbb{R}^{N_\ell}$ and $z_\ell \in \mathbb{R}^{N_\ell}$ to be the values at each layer after and before activation, meaning $a_\ell = \rho(A_\ell a_{\ell-1} + b_\ell)$ and $z_\ell = A_\ell a_{\ell-1}+b_\ell$. Next, consider the quantities
\begin{align}
    \delta_L &= \frac{\partial \mathcal{L}}{\partial A_L} \odot \rho'(z_L) \label{e-deltaL}, \\
    \delta_{\ell} &= ((A_{\ell+1})^T\delta_{\ell+1}) \odot \rho'(z_{\ell}), \text{ } \forall \ell \in \{1, \ldots,L-1\}, \label{e-deltal}
\end{align}
where $\odot$ is the Hadamard product. Note that these values can be efficiently computed via dynamic programming. Moreover, $\forall \ell \in \{1,\ldots,L-1\}$,
\begin{equation}
    \frac{\partial \mathcal{L}}{\partial b_\ell^j} = \delta_\ell^j \text{ and } \frac{\partial \mathcal{L}}{\partial A_\ell^{ij}} = a_{\ell-1}^j \delta_\ell^i, \label{e-grads}
\end{equation}
where superscripts correspond to the index. Using Equations \eqref{e-deltaL}, \eqref{e-deltal}, and \eqref{e-grads}, we now have an efficient way to compute gradients of the loss function with respect to the neural network parameters.

Despite the universal approximation theorem, the neural network defined in Equation~\eqref{d-NN} is not always the best choice of model for certain regression or classification tasks. In practice, one is often able to obtain better results using a neural network with more exotic architectures. For example, one may consider convolutional neural networks (CNNs) or recurrent neural networks (RNNs). In CNNs, the key idea is to replace some of the $W_\ell$ in Equation~\eqref{d-NN} with a convolutional layer $C_\ell$. In a 1D convolutional layer, we pad our input $x \in \mathbb{R}^n$ with $P$ zeros on each side. Then for a layer with stride $S \in \mathbb{N}$ and $F \in \mathbb{N}$ filters with kernel size $K \in \mathbb{N}$, the output is a matrix
\begin{equation}\label{eq:CNN_layer}
    C(x) = \left(\sum_{k = -\lceil K/2 \rceil}^{\lceil K/2 - 1 \rceil} x_{S(i-1)+1-k} h_{j,k} \right)_{1 \leq i \leq n', 1 \leq j \leq F}
\end{equation}
where $(h_{j,k})_{ -\lceil K/2 \rceil}^{\lceil K/2 - 1 \rceil} \in \mathbb{R}^K$ is the $j$-th filter, which we learn during training, and $n' = (n-K+2P)/S + 1$.\footnote{Throughout this paper, we start counting vector indices at 1.} It is also common to follow $C_\ell$ with pooling layers, which replace a block of $M \in \mathbb{N}$ consecutive values with their average or maximum value. For RNNs, we suppose that our data is a sequence of vectors $(x_t)_{t=1}^T$ with $x_t \in \mathbb{R}^n$. We then initialize a hidden state $h_0$, which we update for each $t \in \{1, \ldots, T\}$ according to some learned function of $x_t$. We describe this in more detail in Equation~\ref{eq:GRU_update}.

\begin{figure}[b]
    \centering
    \begin{tabular}{ |p{1.9cm}|p{2.7cm}|p{2.5cm}|p{1.8cm}|p{1.8cm}|p{1.8cm}|}
    \hline
    \multicolumn{6}{|c|}{Low conductor dataset} \\
    \hline
    Degree ($d$) & Conductor ($N$) & \# $L$-functions&  Rank 0 (\%)&  Rank 1 (\%)&  Rank 2 (\%) \\
    \hline
    2 & $\leq 2^{14}$ & 65126 & 41.807\% & 51.879\%& 6.306\%\\
    \hline
    3 & $\leq 2^{15}$ & 41844 & 42.842\% & 51.780\% & 5.374\%\\
    \hline
    4 & $\leq 2^{15}$ & 28987 & 41.932\% & 52.313\% & 5.751\%\\
    \hline
    \end{tabular}\\
    \caption{Data for small neural network experiments.}
    \label{fig:data_dist}
\end{figure}

\section{Datasets}
For the small models described in section 4.1, we will use elliptic curves from the \cite{lmfdb} and a large database of genus 2 curves provided by Andrew Sutherland.\footnote{This data will be added to the LMFDB, but it is not available as of \today.} For degree 2 $L$-functions, we consider $L(E,s)$ where $E$ is a non-CM elliptic curve with conductor $\leq 2^{14}$. For degree 3 $L$-functions, we use $L(E,s+\frac{1}{2})L(\chi,s)$ where $\chi$ is a primitive quadratic Dirichlet character and $N_EN_\chi \leq 2^{15}$. For degree 4 $L$-functions, we consider $L(C, s)$ where $C$ is a genus 2 curves with Sato-Tate group $\text{Usp}(4)$ and conductor $\leq 2^{15}$.

\begin{figure}[t]
    \centering   \includegraphics[width=\textwidth]{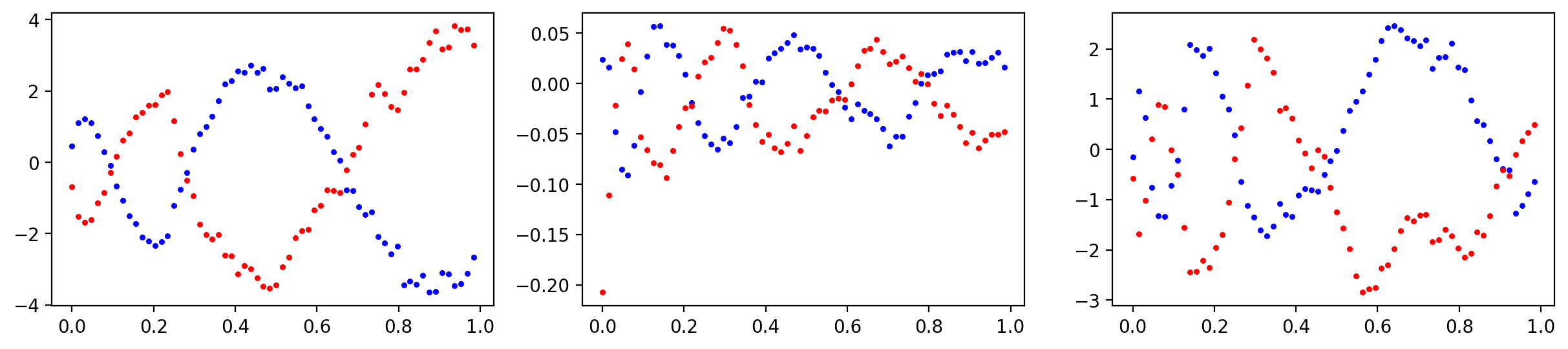}
    \caption{Average entries of $v_{64}(L)$ over all $L$-functions in our low conductor dataset with degree $d$ and root number $w = 1$ (resp. $w = -1$) in \textbf{Blue} (resp. \textbf{Red}) for \textbf{Left:} $d=2$, \textbf{Center:} $d=3$, \textbf{Right:} $d=4$.}
    \label{fig:avg_rn_data}
\end{figure}

\begin{figure}[h]
    \centering    \includegraphics[width=0.5\textwidth]{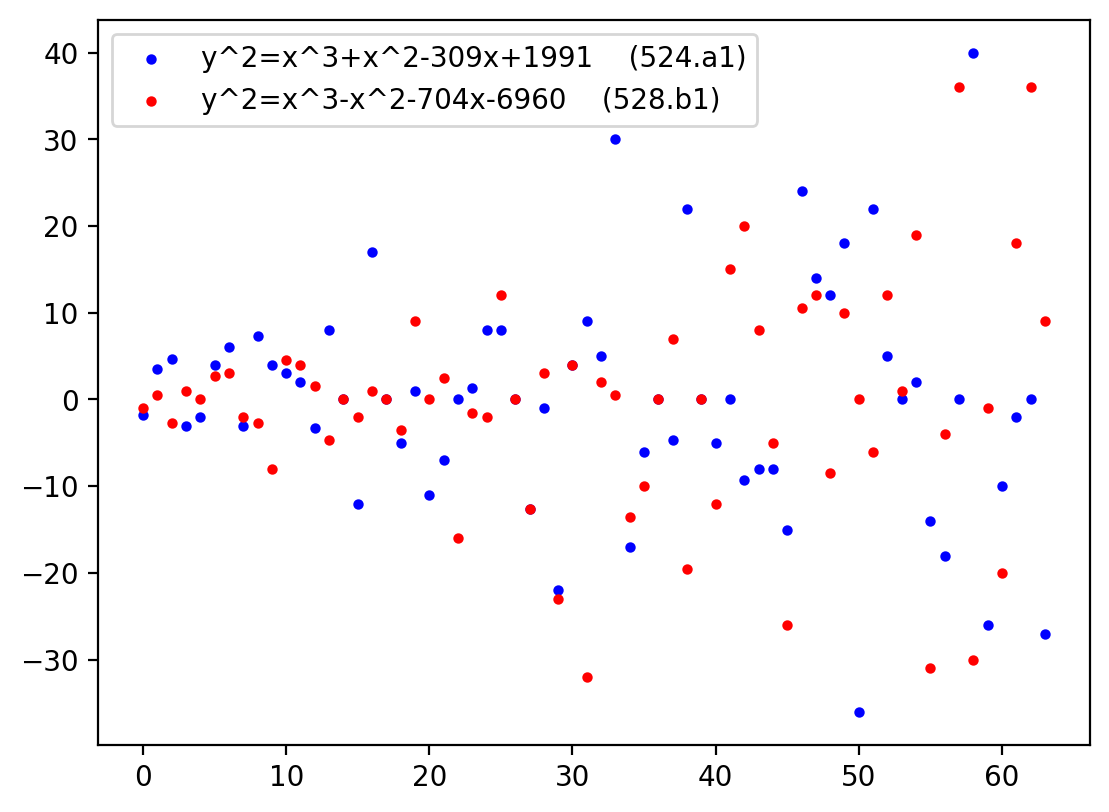}
    \caption{Examples of $v_B(L(E_0, s)$ and $v_B(L(E_1, s)$ where $E_0:y^2=x^3+x^2-309x+1991$ is rank 0 and $E_1:y^2=x^3-704x-6960$ is rank 1.}
    \label{fig:example_vBL}
\end{figure}

Using these $L$-functions, we compute their Dirichlet coefficients $a_p$ and generate a local average inspired by \cite{HLOPS}.\footnote{Such an average was first considered by Johnathan Bober \cite{B23}.} Suppose that $L(s)$ is an $L$-function with Dirichlet coefficients $(a_p)_{p \in \mathcal{P}}$ and conductor $N$. Now fix $B \in \mathbb{N}$ to be some number of bins, and define
\begin{equation}\label{eq:V_B(L)}
    v_{B}(L) = \left(\frac{1}{\#\mathcal{P}_i(B, N)}\sum_{p \in \mathcal{P}_i(B, N)} a_{p} \right)_{i=1}^B \in \mathbb{R}^B,
\end{equation}
where
\begin{equation}
    \mathcal{P}_i(B, N) = \Big\{p \in \mathcal{P} \text{ }\Big| \text{ }\frac{i-1}{B} < \frac{p}{N} \leq \frac{i}{B} \Big\}.
\end{equation}
We use such a construction we two reasons. First, this gives us a standardized input size $B$ which uses all $(a_p)_{p \leq N}$ regardless of $N$. Second, this amplifies the murmuration phenomenon, which is a bias in the $a_p$ that is a function $p/N$ and requires averaging over both $L$-functions and primes to observe. For instance, see Figure~\ref{fig:avg_rn_data} to observe this bias in the dataset described in Figure~\ref{fig:data_dist}. Thus by considering $v_B(L)$ rather than $(a_p)_{p \leq N}$, we hope to build models that generalize to higher conductors more robustly, and that are more likely to pick up on any biases that only appear when averaging over primes. We also provide an example of this construction for two specific elliptic curves in Figure~\ref{fig:example_vBL}.

For the larger models we consider in section 4.2, we generate a random sample of high conductor elliptic curves. Given two positive integers $N_{\text{min}} < N_{\text{max}}$, we find a random elliptic curves $E/\mathbb{Q} : y^2 = x^3 + Ax + B$ such that $A > 0$, $\Delta_E/16 \equiv 1 \bmod 2$, $\mu^2(\Delta_E/16) = 1$, and $N_{\text{min}} \leq N_E = |\Delta_E| \leq N_{\text{max}}$. We note that this is a positive proportion of elliptic curves ordered by height, which follows from \cite{CS20}. To generate such curves, we use the following algorithm:
\begin{enumerate}
    \item Pick a uniform random integer $A \in [1, (N_{\text{max}}/64)^{1/3}]$.
    \item Pick a uniform random odd integer $B \in [\max(1, ((N_\text{min}-64A^3)/432)^{1/2}, ((N_\text{max}-64A^3)/432)^{1/2}]$. If this interval has no odd integers, return to step (1).
    \item If $4A^3 + 27B^2$ is not square-free, return to step (1).
    \item Let $E : y^2 = x^3 + Ax + B$. If $\nu_2(N_E) = 4$, then return $E$. Otherwise, return to step (1).\footnote{$\nu_p$ denotes the $p$-adic valuation.}
\end{enumerate}
Notice that the intervals in steps (1) and (2) ensure that $16(4A^3 + 27B^2) \in [N_{\text{min}}, N_{\text{max}}]$. Moreover, the conditions $B$ odd, $4A^3 + 27B^2$ squarefree, and $\nu_2(N_E) = 4$ ensure that $N_E = |\Delta_E| = 16(4A^3 + 27B^2)$. One can also compute constants $N_0$ and $c$ such that $N_{\text{min}} > N_0$ and $N^{5/6}_{\text{min}} < cN^{5/6}_{\text{max}}$ will provably ensure that the algorithm terminates. In practice, our choices of $N_{\text{min}}$ and $N_{\text{max}}$ will produce more than enough elliptic curves. In particular, by ensuring that all of the curves we find correspond to distinct isogeny classes, we obtain the datasets described in Figures~\ref{fig:mid cond}, \ref{fig:large_cond}. Note that the first row in Figure~\ref{fig:mid cond} corresponds to data from the Stein-Watkins database \cite{SWDB}.

\begin{figure}[h]
    \centering
    \begin{tabular}{ |p{2cm}|p{2.7cm}|p{2.7cm}|p{2cm}|p{2cm}|}
    \hline
    \multicolumn{5}{|c|}{Medium conductor elliptic curve dataset} \\
    \hline
    $N_{\text{min}}$ & $N_{\text{max}}$ & \# $L$-functions& \# $a_p(E)$ &  $w_E = 1$ (\%)\\
    \hline
    $10^{8} - 10^5$ & $10^8$ & 85860&  1229 & 49.886\% \\
    \hline
    $10^{9}$ &  $1.001 \cdot 10^{9}$ & 75200 & 1229 & 49.672\% \\
    \hline
    $10^{10}$ & $1.001 \cdot 10^{10}$ & 93520 & 1229 & 50.029\% \\
    \hline
    \end{tabular}\\
    \caption{Medium conductor data for large neural network experiments.}
    \label{fig:mid cond}
\end{figure}

The goal of the medium conductor dataset will be to better understand how the performance of large CNN models scales with the conductor and $a_p(E)$ count. The specific values in this table were chosen only to ensure we had enough data to see a convincing trend in validation accuracy.

\begin{figure}[h]
    \centering
    \begin{tabular}{ |p{2cm}|p{2.7cm}|p{2.7cm}|p{2cm}|p{2cm}|}
    \hline
    \multicolumn{5}{|c|}{High conductor elliptic curve dataset} \\
    \hline
    $N_{\text{min}}$ & $N_{\text{max}}$ & \# $L$-functions& \# $a_p(E)$ &  $w_E = 1$ (\%)\\
    \hline
    $10^{10}$ & $\approx 1.047997 \cdot 10^{10}$ & 24334&  1457 & 49.490\% \\
    \hline
    $10^{11}$ & $\approx 1.009138 \cdot 10^{11}$ & 31250 & 1821 & 49.914\% \\
    \hline
    $10^{12}$ & $\approx 1.001697 \cdot 10^{12}$ & 39366 & 2227 & 50.196\% \\
    \hline
    $10^{13}$ & $\approx 1.000311 \cdot 10^{13}$ & 48448 & 2707& 50.076\% \\
    \hline
    $10^{14}$ & $\approx 1.000061 \cdot 10^{14}$ & 65536 & 3512& 50.356\% \\
    \hline
    $10^{15}$ & $\approx 1.000012 \cdot 10^{15}$ & 78608 & 4137& 50.164\% \\
    \hline
    \end{tabular}\\
    \caption{High conductor data for large neural network experiments.}
    \label{fig:large_cond}
\end{figure}

The goal of the high conductor dataset is to see if a poly-logarithmic number of $L$-functions and $a_p(E)$ is sufficient to learn root numbers as $N_{\text{min}} \to \infty$. We choose these values to respect the following asymptotic bounds. We use $N_{\text{max}} = N_{\text{min}} + cN_{\text{min}}^{1/6}$ where $c = O(\log(N_{\text{min}})^3)$ so that, assuming the conjecture that $\#\{E/\mathbb{Q} : N_E \leq X\} \sim \alpha X^{5/6}$ for some $\alpha > 0$, we ensure that there is at least $O(\log(N_\text{min})^3)$ curves in our intervals. We then sample $O(\log(N_{\text{min}})^3)$ curves and compute $a_p(E)$ for $p$ up to approximately $\log(N_\text{min})^{3}$.

\section{Models}
\subsection{Small Neural Networks}
In this section, we describe the interpretable models used to classify the data in our low conductor dataset. First, we describe simple CNN models, both in the language of TensorFlow and as explicit mathematical functions. In these models, we will have a 1D convolutional layer with $F \in \{1,2,3\}$ filters and a kernel size of $K = B = 64$. This is followed by a batch normalization layer, which learns to rescale the data to have a mean near 0 and a standard deviation near 1. Then we have a dropout layer, and a reshape layer which just identifies a matrix in $\mathbb{R}^{1 \times F}$ with a vector in $\mathbb{R}^F$. This is followed by a feed-forward dense layer which outputs our prediction, $\hat{w} \in \mathbb{R}^2$ in the case of root numbers and $\hat{r} \in \mathbb{R}^4$ in the case ranks. Note that both the convolutional and dense layers include a bias. See Figure \ref{fig:model_layers_CNN} for an example TensorFlow summary.
\begin{figure}[h]
\begin{center}
\begin{lstlisting}
 Layer                           Output Shape                  Param #   
==============================================================================
 Input                           (64, 1)                       0                                                                                                  
 Conv1D                          (1, 3)                        195                                                           
 Batch Normalization             (1, 3)                        12                                                                                    
 Dropout                         (1, 3)                        0 
 Reshape                         (3)                           0
 Dense (Relu)                    (2)                           8                                        
==============================================================================
Total params: 215
Trainable params: 209
Non-trainable params: 6
\end{lstlisting}    
    \caption{TensorFlow layers for CNN predicting root number with 3 filters.}
    \label{fig:model_layers_CNN}
\end{center}
\end{figure}

Since our kernel size is equal to the input size, this convolutional layer acts more like an ensemble of $F$ dense layers than a true convolutional layer. In this case, we can write out the entire model for predicting root numbers (resp. ranks) explicitly as:
\begin{equation}
    f(v_B(L)) = \sigma\left(A  \left(\frac{\gamma_i \left(\rho( w_i \cdot v_B(L) + b_i) - \mu_i\right)}{\sqrt{\sigma_i^2 + \epsilon} + \beta_i}\right)_{i=1}^F  + b\right),
\end{equation}
where $v_B(L)$ is defined in Equation~\eqref{eq:V_B(L)}, $A \in \mathbb{R}^{F \times 2}$ (resp. $A \in \mathbb{R}^{F \times 4}$), $w_i \in \mathbb{R}^B$, $b_i, \beta_i, \gamma_i, \mu_i \in \mathbb{R}$, $b \in \mathbb{R}^2$ (resp. $b \in \mathbb{R}^4$) $\sigma_i^2 \in \mathbb{R}_{\geq 0}$, $\epsilon = 0.001$, $\rho(x) = \max(0, x)$, and $\sigma$ is a softmax. This gives us a vector of size 2 (resp. 4), whose $j$-th entry we interpret as the probability that $L$ has root number $(-1)^{j+1}$ (resp. $L$ has rank $j-1$). We also consider a slight variation of this model where we use the truncated vectors 
\begin{equation}
    v'_{B}(L) = \left(\frac{1}{|\mathcal{P}_i(B, N)|}\sum_{p \in \mathcal{P}_i(B, N)} a_{p} \right)_{i=5}^B \in \mathbb{R}^{B-4},
\end{equation}
discarding $(a_p)_{p \leq \frac{4N}{B}}$. In the following figures and discussion, we refer to models that use $v_B(L)$ as the full models and those that use $v'_B(L)$ as the partial models.

We also consider an even smaller RNN model. We will be using a gated recurrent unit (GRU), which is introduced in \cite{C14}. This is an RNN where we initialize the hidden state to $h_0 = 0 \in \mathbb{R}^3$, and for $t \in \{1, \ldots, T\}$, we update this state inductively by the rules
\begin{align}
\begin{split}\label{eq:GRU_update}
    z_t &=\sigma(W_z x_t + U_z h_{t-1} + b_z), \\
    r_t &=\sigma(W_r x_t + U_r h_{t-1} + b_r), \\
    \hat{h}_t &=\tanh(W_h x_t + U_t(r_t \odot h_{t-1}) + b_h),\\
    h_t &= (1-z_t) \odot h_{t-1} + z_t \odot \hat{h}_{t},
\end{split}
\end{align}
where $\sigma: \mathbb{R} \to \mathbb{R}$ is a sigmoid function, and our weights are the entries of the following matrices and vectors: $W_z, W_r, W_h \in \mathbb{R}^{3 \times 1}$, $U_z, U_r, U_h \in \mathbb{R}^{3 \times 3}$, and $b_z, b_r, b_h \in \mathbb{R}^3$. The GRU then outputs $h_T$, which we input to a feed-forward dense layer. Note that in this case, we treat $v_B(L)$ as a sequence of $B$ real numbers.

\begin{figure}[ht]
\begin{center}
\begin{lstlisting}
 Layer                           Output Shape                  Param #   
==============================================================================
 Input                           (64, 1)                       0           
 GRU                             (3)                           45                     
 Dropout                         (3)                           0 
 Dense (Relu)                    (2)                           8                                        
==============================================================================
Total params: 53
Trainable params: 53
Non-trainable params: 0
\end{lstlisting}    
    \caption{TensorFlow layers for RNN predicting root numbers with 3 units.}
    \label{fig:model_layers}
\end{center}
\end{figure}

\begin{figure}[p]
    \centering
    \includegraphics[width=\textwidth]{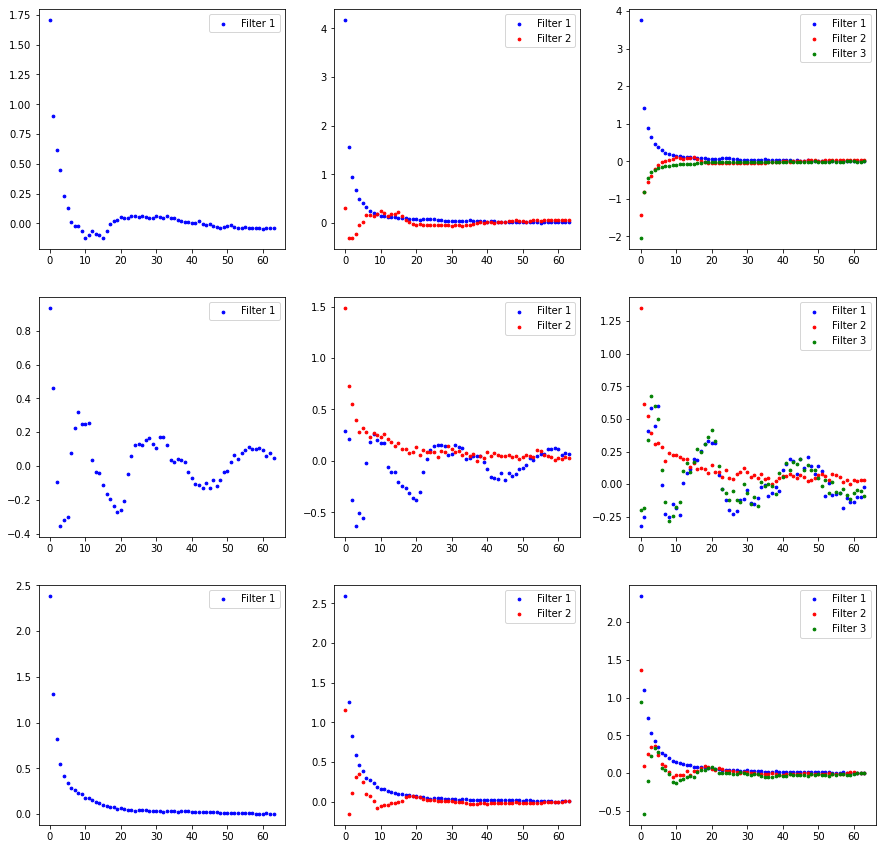}
    \caption{Learned filters in full CNN models for predicting root numbers corresponding to degree $d = 2, 3, 4$ (top, middle, bottom) with $F = 1, 2, 3$ (left, middle, right).}
    \begin{tabular}{ |p{1.8cm}|p{2.4cm}|p{2.5cm}|p{4cm}|p{4cm}|}
    \hline
    \multicolumn{5}{|c|}{Accuracy of full CNN root number classifiers} \\
    \hline
    Degree ($d$)& \#Filters ($F$) & Test Accuracy & Generalization Accuracy 1 & Generalization Accuracy 2 \\
    \hline
    2 & 1 & 0.9368 & 0.9192 & 0.9086\\
    \hline
    2 & 2 & 0.9763 & 0.9538 & 0.9369\\
    \hline
    2 & 3 & 0.9906 & 0.9715 & 0.9497\\
    \hline
    3 & 1 & 0.9541 & 0.9385 & 0.9002\\
    \hline
    3 & 2 & 0.9818 & 0.9369 & 0.7062\\
    \hline
    3 & 3 & 0.9827 & 0.9416 & 0.6973\\
    \hline
    4 & 1 & 0.9327 & 0.8762 & 0.8093\\
    \hline
    4 & 2 & 0.9717 & 0.9398 & 0.9020\\
    \hline
    4 & 3 & 0.9721 & 0.9425 & 0.9052\\
    \hline
    \end{tabular}
    \label{fig:rn_full}
\end{figure}

\begin{figure}[p]
    \centering
    \includegraphics[width=\textwidth]{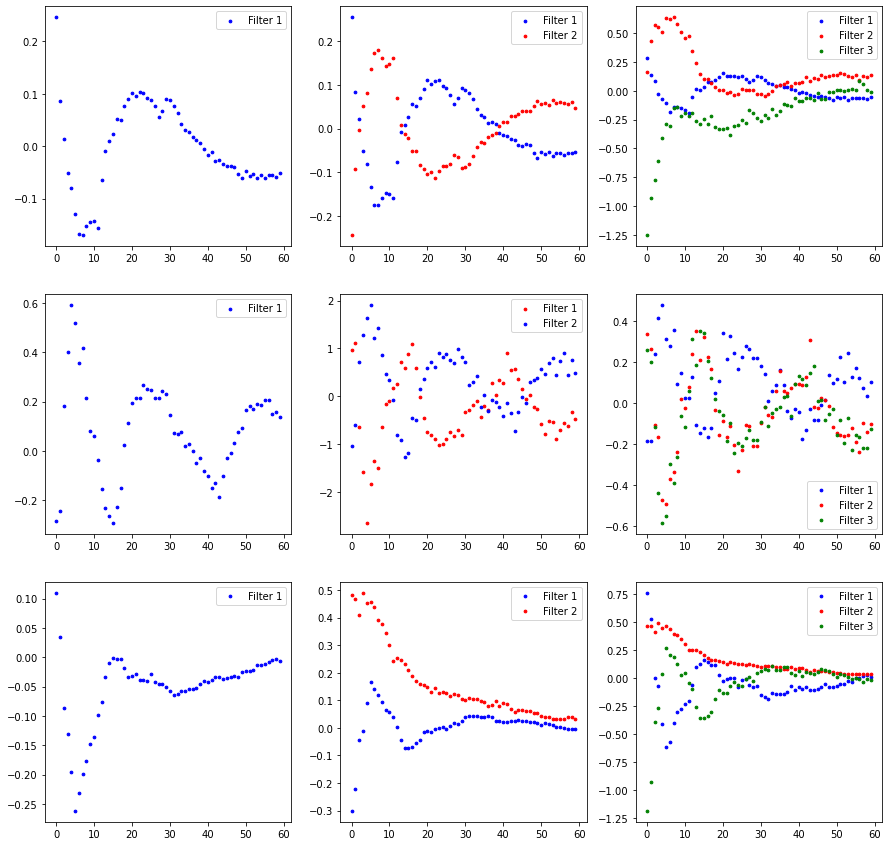}
    \caption{Learned filters in partial CNN models for predicting root numbers corresponding to degree $d = 2, 3, 4$ (top, middle, bottom) with $F = 1, 2, 3$ (left, middle, right).}
    \begin{tabular}{ |p{1.8cm}|p{2.4cm}|p{2.5cm}|p{4cm}|p{4cm}|}
    \hline
    \multicolumn{5}{|c|}{Accuracy of partial CNN root number classifiers} \\
    \hline
    Degree ($d$)& \#Filters ($F$) & Test Accuracy & Generalization Accuracy 1 & Generalization Accuracy 2 \\
    \hline
    2 & 1 & 0.8761 & 0.8662 & 0.8586 \\
    \hline
    2 & 2 & 0.8824 & 0.8733 & 0.8655 \\
    \hline
    2 & 3 & 0.8856 & 0.8730 & 0.8650 \\
    \hline
    3 & 1 & 0.9165 & 0.9085 & 0.8872 \\
    \hline
    3 & 2 & 0.9176 & 0.9043 & 0.8895 \\
    \hline
    3 & 3 & 0.9178 & 0.9032 & 0.8864 \\
    \hline
    4 & 1 & 0.8423 & 0.7903 & 0.7530 \\
    \hline
    4 & 2 & 0.8622 & 0.8118 & 0.7737 \\
    \hline
    4 & 3 & 0.8584 & 0.8144 & 0.7773 \\
    \hline
    \end{tabular}
    \label{fig:rn_partial}
\end{figure}

\begin{figure}[p]
    \centering
    \includegraphics[width=\textwidth]{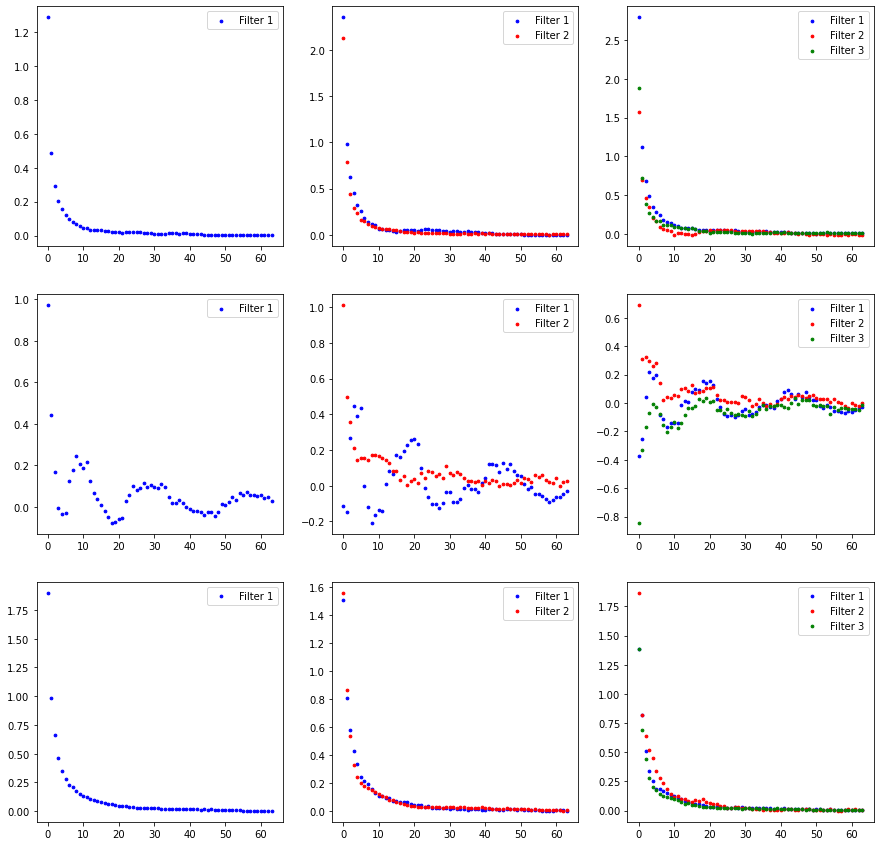}
    \caption{Learned filters in full CNN models for predicting ranks corresponding to degree $d = 2, 3, 4$ (top, middle, bottom) with $F = 1, 2, 3$ (left, middle, right).}
    \begin{tabular}{ |p{1.8cm}|p{2.4cm}|p{2.5cm}|p{4cm}|p{4cm}|}
    \hline
    \multicolumn{5}{|c|}{Accuracy of full CNN rank classifiers} \\
    \hline
    Degree ($d$)& \#Filters ($F$) & Test Accuracy & Generalization Accuracy 1 & Generalization Accuracy 2 \\
    \hline
    2 & 1 & 0.9765 & 0.9279 & 0.8753\\
    \hline
    2 & 2 & 0.9900 & 0.9712 & 0.9420\\
    \hline
    2 & 3 & 0.9900 & 0.9726 & 0.9488\\
    \hline
    3 & 1 & 0.9406 & 0.8857 & 0.7762\\
    \hline
    3 & 2 & 0.9762 & 0.9236 & 0.6728\\
    \hline
    3 & 3 & 0.9804 & 0.9151 & 0.7601\\
    \hline
    4 & 1 & 0.9880 & 0.9518 & 0.9026\\
    \hline
    4 & 2 & 0.9955 & 0.9792 & 0.9570\\
    \hline
    4 & 3 & 0.9952 & 0.9783 & 0.9552\\
    \hline
    \end{tabular}\\
    \label{fig:rank_full}
\end{figure}

\begin{figure}[p]
    \centering
    \includegraphics[width=\textwidth]{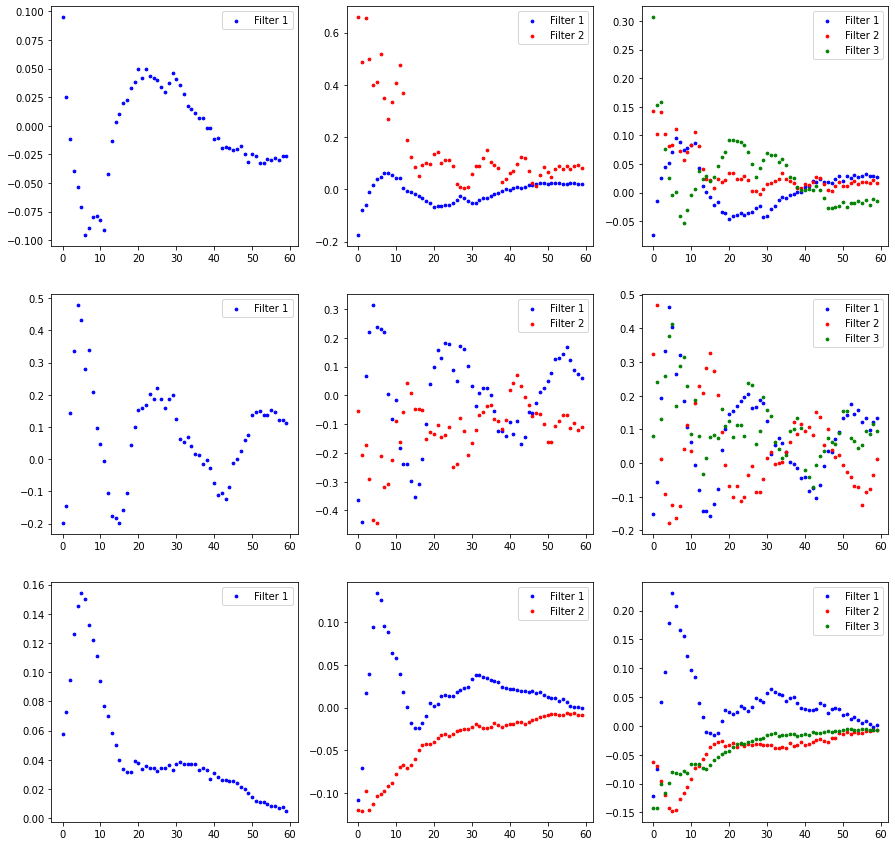}
    \caption{Learned filters in partial CNN models for predicting ranks corresponding to degree $d = 2, 3, 4$ (top, middle, bottom) with $F = 1, 2, 3$ (left, middle, right).}
    \begin{tabular}{ |p{1.8cm}|p{2.4cm}|p{2.5cm}|p{4cm}|p{4cm}|}
    \hline
    \multicolumn{5}{|c|}{Accuracy of partial CNN rank classifiers} \\
    \hline
    Degree ($d$)& \#Filters ($F$) & Test Accuracy & Generalization Accuracy 1 & Generalization Accuracy 2 \\
    \hline
    2 & 1 & 0.8197 & 0.6948 & 0.6564 \\
    \hline
    2 & 2 & 0.8550 & 0.7751 & 0.7427 \\
    \hline
    2 & 3 & 0.8540 & 0.7774 & 0.7440 \\
    \hline
    3 & 1 & 0.8678 & 0.7609 & 0.7150 \\
    \hline
    3 & 2 & 0.8736 & 0.7645 & 0.7194 \\
    \hline
    3 & 3 & 0.8749 & 0.7674 & 0.7183 \\
    \hline
    4 & 1 & 0.8309 & 0.7299 & 0.6469 \\
    \hline
    4 & 2 & 0.8596 & 0.7767 & 0.7036 \\
    \hline
    4 & 3 & 0.8601 & 0.7784 & 0.7047 \\
    \hline
    \end{tabular}\\
    \label{fig:rank_partial}
\end{figure}

\begin{figure}[t]
    \centering
    \includegraphics[width=\textwidth]{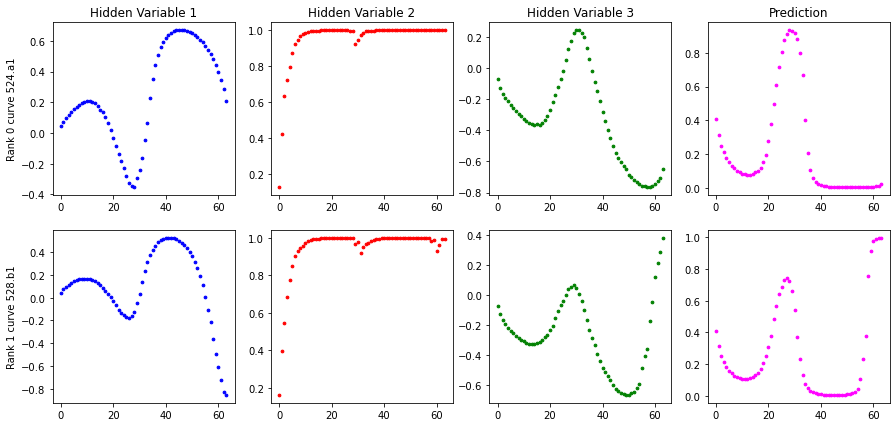}
    \caption{Hidden state and prediction in GRU as a function of $t \in \{1, \ldots, B\}$ for elliptic curves $E_0:y^2=x^3+x^2-309x+1991$ (rank 0) and $E_1:y^2=x^3-704x-6960$ (rank 1).}
    \begin{tabular}{|p{1.8cm}|p{2.5cm}|p{4cm}|p{4cm}|}
    \hline
    \multicolumn{4}{|c|}{Accuracy of full RNN root number classifiers} \\
    \hline
    Degree ($d$)& Test Accuracy & Generalization Accuracy 1 & Generalization Accuracy 2 \\
    \hline
    2 &  0.9807 & 0.9439 & 0.9047 \\
    \hline
    3 & 0.9822 & 0.9515 & 0.8314 \\
    \hline
    4 & 0.9919 & 0.9774 & 0.9557 \\
    \hline
    \end{tabular}
    \newline
    \begin{tabular}{|p{1.8cm}|p{2.5cm}|p{4cm}|p{4cm}|}
    \hline
    \multicolumn{4}{|c|}{Accuracy of partial RNN root number classifiers} \\
    \hline
    Degree ($d$)& Test Accuracy & Generalization Accuracy 1 & Generalization Accuracy 2 \\
    \hline
    2 &  0.8740 & 0.8624 & 0.8556 \\
    \hline
    3 & 0.9085 & 0.9022 & 0.8907 \\
    \hline
    4 & 0.8439 & 0.7908 & 0.7543 \\
    \hline
    \end{tabular}
    \newline
    \begin{tabular}{|p{1.8cm}|p{2.5cm}|p{4cm}|p{4cm}|}
    \hline
    \multicolumn{4}{|c|}{Accuracy of full RNN rank classifiers} \\
    \hline
    Degree ($d$)& Test Accuracy & Generalization Accuracy 1 & Generalization Accuracy 2 \\
    \hline
    2 &  0.9863 & 0.9588 & 0.9180 \\
    \hline
    3 & 0.9915 & 0.9219 & 0.7764 \\
    \hline
    4 & 0.9912 & 0.9739 & 0.9469 \\
    \hline
    \end{tabular}
    \newline 
    \begin{tabular}{|p{1.8cm}|p{2.5cm}|p{4cm}|p{4cm}|}
    \hline
    \multicolumn{4}{|c|}{Accuracy of partial RNN rank classifiers} \\
    \hline
    Degree ($d$)& Test Accuracy & Generalization Accuracy 1 & Generalization Accuracy 2 \\
    \hline
    2 &  0.8165 & 0.6956 & 0.6586 \\
    \hline
    3 & 0.8679 & 0.7637 & 0.7222 \\
    \hline
    4 & 0.8203 & 0.7188 & 0.6427 \\
    \hline
    \end{tabular}
    \newline
    \label{fig:RNN_model}
\end{figure}

For both the CNNs and RNNs, we train the models on a random subset of 80\% of the data in Figure~\ref{fig:data_dist}. For CNNs, we use the Adam optimizer with a cyclic learning rate \cite{Sm15} for 1,000 epochs and a dropout of 0.1. For RNNs, we use Adam for 10,000 epochs and a dropout rate of 0.2. We then report on the accuracy of these models on the remaining 20\% of the data, as well as two generalization accuracies. If our data in Figure~\ref{fig:data_dist} for degree $d$ has conductor up to $2^n$, we call generalization accuracy 1 (resp. 2) the accuracy of the model on $L$-functions with conductor in $(2^n, 2^{n+1}]$ (resp. $(2^{n+1}, 2^{n+2}]$). We plot the learned filters and report on the accuracy of the CNN models in Figures~\ref{fig:rn_full}, \ref{fig:rn_partial}, \ref{fig:rank_full}, \ref{fig:rank_partial}. We plot the hidden states in the GRU obtained from applying the model to the curves in Figure~\ref{fig:example_vBL} and report on the accuracy of the RNN models in Figure~\ref{fig:RNN_model}.

\subsection{Large Neural Networks}
In this section, we describe the models used to predict root numbers of our medium and high conductor datasets of elliptic curves. For both datasets, we use large CNNs with architectures inspired by \cite{KV22}. In particular, we will have $n_I \in \mathbb{Z}_{\geq 0}$ input convolutional blocks, $n_R \in \mathbb{Z}_{\geq 0}$ reducing convolutional blocks, and then $n_O \in \mathbb{Z}_{\geq 0}$ output convolutional blocks. In each convolutional block, we have a 1D convolutional layer, described in Equation~\eqref{eq:CNN_layer}, followed by an average pooling layer, a batch normalization layer, and a dropout layer. Note that all convolutional layers have $F$ filters with a kernel size $K$ and activation function $\rho(x) = \max(0, x)$. The input and output convolutional blocks have a stride $S = 1$, and the reducing convolutional blocks have a stride of $S=2$. After the convolutional blocks, we have one feed-forward layer with softmax activation. Every convolutional and dense layer includes a bias.

\begin{figure}[ht]
    \centering
    \includegraphics[width=\textwidth]{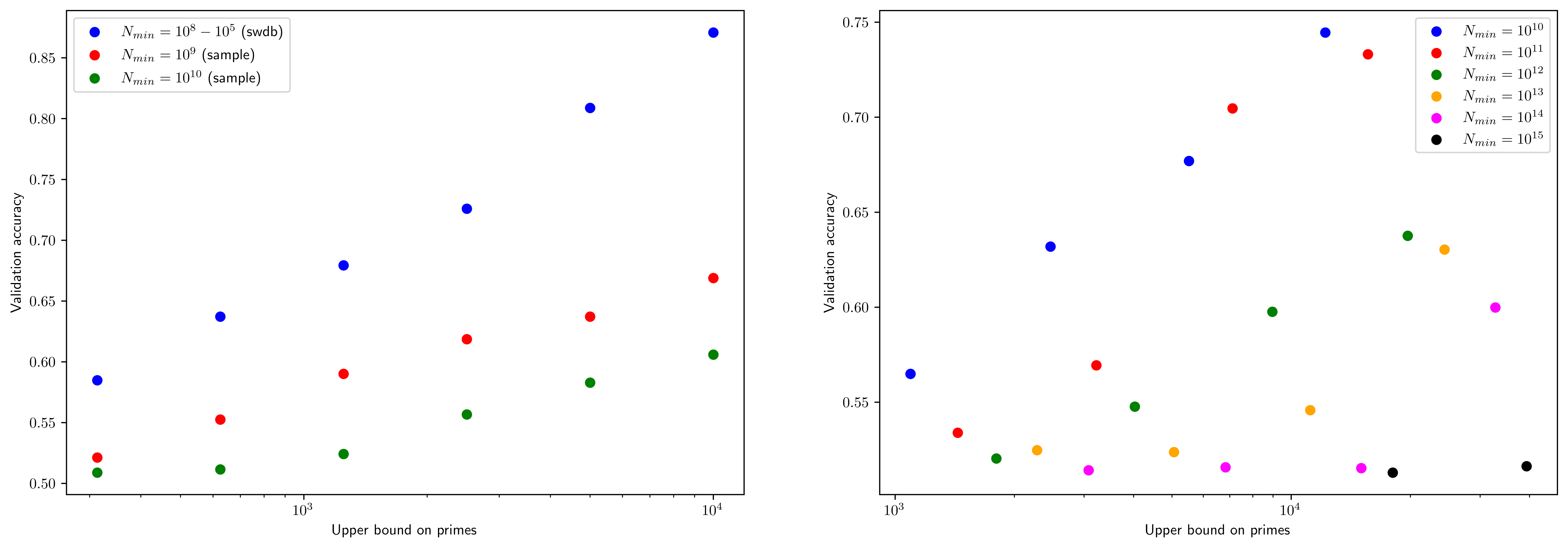}
    \caption{Validation accuracy of large CNN models as we vary the number of $a_p(E)$ and $N_{\text{min}}$.}
    \label{fig:large_CNN}
\end{figure}

Using such a model, we performed the following classification problem. For each conductor interval in Figure~\ref{fig:mid cond}, we attempt to predict the root number of an elliptic curve $E$ using $a_p(E)$ with $p$ up to $10^4/2^i$ for $i \in \{0, \ldots, 5\}$. For these experiments, we had a fixed pool size of $4$ and no padding. All other hyper-parameters were manually tuned for each model. Likewise, for each conductor interval in Figure~\ref{fig:large_cond}, we attempt to predict the root number of an elliptic curve $E$ using $a_p(E)$ with $p$ up to approximately $\log(N_{\text{min}})^3$, as well as using only the first $1/2$, $1/4$, and $1/8$ of the $a_p(E)$.\footnote{See Figure~\ref{fig:large_cond} for the exact number.} For these models, we include padding and we set $n_I = 0$ since this was found to be optimal in \cite{KV22}. In an attempt to systematically maximize our model performance, we perform a hyperparameter search using Bayesian optimization for each of these classification problems. We search for the filter count $F \in \{8, \ldots, 32\}$, kernel size $K \in \{8, \ldots, 64\}$, the pooling size $M \in \{1, \ldots, 8\}$, the dropout rate $q \in \{0.1, 0.2, 0.3\}$, the reducing block count $n_R \in \{0,\ldots, \lfloor \log_2(\log(N_E)^3/\log(\log(N_E)^3)) \rfloor\}$, the output block count $n_O \in \{0, \ldots, 6\}$, and the learning rate $l_r \in [5\cdot 10^{-4}, 5\cdot10^{-2}]$ which maximize the validation accuracy of the CNN. For our Bayesian optimization, we do a total of 128 iterations where the first 32 are random initial trials. The models we obtain from Bayesian optimization typically have a parameter count in $\approx [10^5, 10^6]$. We plot the validation accuracies of the resulting models in Figure~\ref{fig:large_CNN}.

\newpage

\section{Results}
In this section, we provide an analysis of our results, give relevant heuristics, and in the spirit of \cite{KV22}, leave questions for further research. We begin by observing two clear influences in the filters learned by the small CNN models. The first influence is coming from Mestre-Nagao heuristics. We note that a standard Mestre-Nagao sum, which was analyzed in detail in \cite{KM22}, is given by
\begin{equation}\label{eq:Mestre-Nagao}
    \lim_{X \to \infty} \frac{1}{\log X} \sum_{\substack{p \nmid N_E \\ p \leq X}} \frac{a_p(E)\log p }{p} = -r_E + \frac{1}{2}.
\end{equation}
When properly rescaled, this $\log p/p$ weighting matches closely with the filters learned in the degree 2 and 4 models in Figure~\ref{fig:rank_full}. The second influence comes from murmurations, which can be seen in the filters learned in Figure~\ref{fig:rn_partial}. From the full models, we see that the Mestre-Nagao weighting is the dominant weighting in degrees 2 and 4, but the murmuration weighting appears to be best for degree 3. This is quite different from the partial models where we see the murmuration weightings in most filters. The decrease in accuracy observed in the partial models suggests that the Mestre-Nagao weighting at small primes is crucial for high-accuracy predictions, particularly for degrees 2 and 4. This is consistent with the observation that partial models seem to perform best in degree 3, where the murmuration weighting plays a more significant role in the full models.

We also observe that the murmuration weighting is more relevant for predicting root numbers than rank. For instance, the top left plot in Figure~\ref{fig:rn_full} appears to be a combination of the two weightings, which become separated in the top middle plot where we allow the model to learn two filters. Compare this with the degree 2 plots in Figure~\ref{fig:rank_full} where the murmuration weighting appears to be much more subtle. What is perhaps unexpected is that the CNNs often obtain better accuracy in predicting ranks than root numbers, even though one can easily recover the root number from the rank. While this suggests predicting ranks may be a better approach to predicting root numbers, the Goldfeld-Katz-Sarnak conjecture that 50\% of elliptic curves have rank 0 and 50\% of elliptic curves have rank 1 when ordered by conductor prevents this from being a relevant distinction for large enough conductors.

Next, we compare our RNN and CNN models. As seen in Figure~\ref{fig:RNN_model}, the value of the hidden variables may drastically change when the model sees the $a_p$ for $p$ near the top of our range. While this may suggest that the RNN models are less sensitive than CNN models to the $a_p$ at small primes, we see an almost identical decrease in accuracy when considering partial RNN models. These two models are also fundamentally different in that the CNN weights vary with $p/N_E$, meanwhile, the RNN weights are independent of $p/N_E$. The similarity in accuracy between these two approaches seems to suggest there is some upper bound to the accuracy of statistical learning methods. We get a better understanding of this upper bound from our medium and high conductor elliptic curve experiments, which we analyze next.

The left plot in Figure~\ref{fig:large_CNN} suggests that if we use $(a_p)_{p \leq X}$, then our validation accuracy for a fixed $N_{\text{min}}$ is
\begin{equation}\label{eq:Acc}
    \text{Acc}(X) \approx \begin{cases}
        0.5 & \text{if } X < X_{\text{min}} \\
        0.5 + \frac{\log X - \log X_{\text{min}}}{2(\log X_{\text{max}} - \log X_\text{min})} & \text{if } X_{\text{min}} \leq X \leq X_{\text{max}} \\ 
        1 & \text{if } X > X_{\text{max}}
    \end{cases}
\end{equation}
where $X_{\text{min}}$ and $X_{\text{max}}$ depend on our conductor range. This equation also seems to hold for most of our high conductor models, although it is less clear in the $N_{\text{min}} = 10^{11}$ dataset. In particular, more data is needed to be confident that the $\log X$ growth persists in general. Nevertheless, we see a relatively steady increase in accuracy as a function of $\log X$. From the right plot in Figure~\ref{fig:large_CNN}, we observe that $X = O(\log(N_{\text{min}})^3)$ appears insufficient as $N_{\text{min}} \to \infty$. To understand these observations theoretically, we compare this to Mestre-Nagao type heuristics and a probabilistic model. 

One Mestre-Nagao sum considered in \cite{KV22}, which is based on the explicit formula and first worked out in \cite{B13}, is given by
\begin{align}\label{eq:Bober_sum}
\begin{split}
    S(Y) =  \frac{\log N_E}{2\pi Y} &- \frac{\log 2 \pi}{\pi Y} - \frac{1}{\pi Y} \sum_{p \leq \exp(2\pi Y)} \sum_{k=1}^{\lceil 2\pi Y/\log p \rceil} \frac{c_{p^k} \log p }{p^{k/2}}\left(1 - \frac{k \log p}{2 \pi Y}\right) \\ &+ \frac{1}{\pi} \Re\left\{\int_{-\infty}^\infty \frac{\Gamma'(1+it)}{\Gamma(1+it)}\left(\frac{\sin(\pi Y t)}{\pi Y t} \right)^2dt\right\},
\end{split}
\end{align}
where
\begin{equation}
 c_{n} = \begin{cases}
     \alpha_p^k + \beta_p^k & \text{ if } n = p^k \text{ and } p\nmid N_E, \\
     a_p^k & \text{ if } n = p^k \text{ and } p \mid N_E,\\
     0 & \text{ otherwise},
 \end{cases}
\end{equation}
and $\alpha_p, \beta_p$ are eigenvalues of the Frobenius morphism at $p$. Bober proved that if $L(E,s)$ satisfies the Riemann hypothesis, than $\lim_{Y \to \infty} S(Y) = \text{ord}_{s=1} L(E,s)$. Based on the very crude bound that
\begin{equation}
    \frac{\Gamma'(1+it)}{\Gamma(1+it)} = O(t^{1/2}),
\end{equation}
it is easy to see that 
\begin{equation}
    \frac{1}{\pi} \Re\left\{\int_{-\infty}^\infty \frac{\Gamma'(1+it)}{\Gamma(1+it)}\left(\frac{\sin(\pi Y t)}{\pi Y t} \right)^2dt\right\} \ll \frac{1}{Y} \int_{-\infty}^\infty \frac{\sin(\pi Y t)^2}{\pi Y |t|^{3/2}}dt = 2Y^{-3/2}.
\end{equation}
If we set $\log X = 2 \pi Y$, then we have
\begin{equation}
    S(X) = \frac{\log N_E}{\log X} - \frac{2}{\log X} \sum_{p \leq X} \sum_{k=1}^{\lceil \log X/\log p \rceil} \frac{c_{p^k} \log p }{p^{k/2}}\left(1 - \frac{k \log p}{\log X}\right) + O\left(\frac{1}{\log X}\right), 
\end{equation}
as $X \to \infty$. We conclude that 
\begin{equation}\label{eq:heuristic1}
    2 \sum_{p \leq X} \sum_{k=1}^{\lceil \log X/\log p \rceil} \frac{c_{p^k} \log p }{p^{k/2}}\left(1 - \frac{k \log p}{\log X}\right) = - r_E \log X + \log N_E + O(1).
\end{equation}
Thus we need $\log X \gg \log N_E$ to start detecting $r_E$. This implies that need $X \gg N_E^c$ for some $c > 0$, and that $X = O(\log(N_E)^\alpha)$ will be insufficient for any $\alpha > 0$. We also note that the truncated Perron's formula \cite[Proposition 4]{KM22} and zero counts for $L(E,s)$ \cite[Theorem 5.5.9]{Sp15} suggest that error terms in more standard Mestre-Nageo sums (e.g. Equation~\eqref{eq:Mestre-Nagao}) depend on $\log N_E$, although this is less clear than Equation~\eqref{eq:heuristic1}.

Now we present a probabilistic heuristic, noting that such a model is particularly appropriate in the context of statistical learning.\footnote{This heuristic was first worked out by Andrew Granville, who attributes it to a discussion with James Maynard.} Suppose that $(a_p(E))_{p \in \mathcal{P}}$ are independent random variables with mean $\mu_E = 1/2 - r_E$ and variance $\asymp p$. Then by Bienaymé's identity and the prime number theorem,
\begin{equation}
    \frac{1}{\pi(X)}\sum_{X < p \leq 2X} a_p(E) \text{ has variance $\asymp \log X$, }
\end{equation}
 as $X \to \infty$. Since $\pi(X) \asymp p/\log p$ in this range, this implies
\begin{equation}
    \sum_{X < p \leq 2X} \frac{a_p(E) \log p}{p} \text{ has variance $\asymp \log X$.}
\end{equation}
Setting $X = 2^j$, this sum has variance $\asymp j$. Averaging over $J < j \leq 2J$, we get that
\begin{equation}\label{eq:prob_sum_1}
    \frac{1}{\log X} \sum_{X < p \leq X^2} \frac{a_p(E) \log p}{p} \text{ has variance $\asymp  1$. }
\end{equation}
Since $\log X \asymp \log p$ in this range, this implies
\begin{equation}
    \sum_{X < p \leq X^2} \frac{a_p(E)}{p} \text{ has variance $\asymp 1$.}
\end{equation}
Setting $J = 2^k$ and averaging over $K < k \leq 2K$ gives us a random variable with variance $1/K$. In particular, for $m \asymp 2^k$ large enough, we expect that
\begin{equation}\label{eq:prob_sum_2}
    \frac{1}{\log \log X}\sum_{X < p \leq X^m} \frac{a_p(E)}{p} \text{ has variance $\asymp \frac{1}{\log m}$}.
\end{equation}
Therefore, to see convergence we expect to need an exponentially large prime range. We note that the sums in Equations~\eqref{eq:prob_sum_1} and \eqref{eq:prob_sum_2} are both truncated versions of Mestre-Nagao sums appearing in \cite{S07}. 

This heuristic also suggests that using other ranges of primes will not aid in our search for a low-complexity statistic of $a_p(E)$. Notice that Equation~\eqref{eq:prob_sum_2} implies that
\begin{equation}
    \sum_{Y < p \leq X} \frac{a_p(E)}{p} \Bigg/ \sum_{Y < p \leq X} \frac{1}{p} 
\end{equation}
will only be near the mean when $\log X/ \log Y$ is sufficiently large. Based on \cite{HLOPS}, one might expect that $a_p(E)$ for $p \in [N_E/25, N_E/25 + O(\log(N_E)^{3})]$, which is near where the murmuration density for elliptic curves reaches its first peak, would be a good input to our models. Using this set of primes and performing the same hyperparameter search that was used to generate Figure~\ref{fig:large_CNN} gave us a validation accuracy of 0.5290 on the $N_{\text{min}} = 10^{10}$ dataset, and lower accuracy on larger conductors. Based on these heuristics and our experimental results, it seems unlikely that one will be able to recover a low-complexity statistic for the prediction of root numbers from a straightforward machine learning approach. We conclude with a few ideas for future work.
\begin{enumerate}
    \item We saw that machine learning models behave differently in degree 3. Could one take advantage of degree 3 $L$-functions to obtain higher accuracy predictions of root numbers and ranks of elliptic curves? For instance, can we get better predictions for $w_E$ from the Dirichlet coefficients of $\{L(\chi,s)L(E, s + \frac{1}{2})\}_{\chi \in \mathcal{Q}}$ where $\mathcal{Q}$ is a suitable set of primitive quadratic Dirichlet characters?
    \item While a straightforward application of machine learning was unable to predict root numbers in polynomial time, there may be a more sophisticated approach that yields such an algorithm. Namely, is there a polynomial-time reduction of root number prediction which machine learning can solve in polynomial time?
\end{enumerate}

%%%%%%%%%%%%%%%%%%%%%%%%%
%%%==============================

\end{document}